\documentclass[12pt]{article}

 \usepackage{amsfonts,amssymb,amsthm,amsmath,graphicx,mathrsfs,color}

\usepackage[pdfstartview=FitH,
            CJKbookmarks=true,
            bookmarksnumbered=true,
            bookmarksopen=true,
            colorlinks,
            linkcolor=blue,
            anchorcolor=blue,
            citecolor=blue,
            urlcolor=blue
            ]{hyperref}

\newtheorem{theo}{Theorem}
\newtheorem{lemma}[theo]{Lemma}
\newtheorem{coro}[theo]{Corollary}
\newtheorem{prop}[theo]{Proposition}
\newtheorem{conj}[theo]{Conjecture}
\newtheorem{example}[theo]{Example}
\newtheorem{remark}[theo]{Remark}

\def\Z{\mathbb Z}

\def\N{\mathbb N}
\def\R{\mathbb R}

\begin{document}
\author{K\'aroly J. B\"or\"oczky\footnote{funded by NKFIH grants 132002 and 150613.}, M\'aty\'as Domokos\footnote{funded by NKFIH grants 132002 and 138828.}, Ansgar Freyer\footnote{funded by the Deutsche Forschungsgemeinschaft (DFG, German Research Foundation) grant 539867386.}, \\Christoph Haberl,
 Gergely Harcos\footnote{funded by Hungarian Academy of Sciences "Lend\"ulet" grant "Automorphic Research Group" and NKFIH grant 143876.}, Jin Li\footnote{funded by the National Natural Science Foundation of China grant 12201388 and the Austrian Science Fund (FWF) grant I3027.}}

\title{Exponential valuations on lattice polygons}


\maketitle

\begin{abstract}
We classify translatively exponential and ${\rm GL}(2,\Z)$ covariant valuations on lattice polygons with values in the space of real (complex) measurable functions. A typical example of such valuations is induced by the Laplace transform, but as it turns out there are many more. The argument uses the ergodicity of the linear action of ${\rm SL}(2,\Z)$ on $\R^2$, and some elementary properties of the Fibonacci numbers.
\end{abstract}
\bigskip


\section{Exponential Valuations}

By a {\it polytope} $P\subset\R^n$, we mean the convex hull of finitely many points of $\R^n$. The polytope $P$ is called a {\it segment}
if ${\rm dim}\,P=1$, and $P$ is called a {\it polygon}
if ${\rm dim}\,P=2$.
Let $\mathcal{F}$ be a family of compact convex sets in $\R^n$ such that if $P\cup Q$ is convex for $P,Q\in\mathcal{F}$,
then $P\cup Q\in\mathcal{F}$ and $P\cap Q\in\mathcal{F}$. Examples of such families are the
   family $\mathcal{K}^n$ of all convex compact sets in $\R^n$, the family
$\mathcal{P}^n$ of all polytopes in $\R^n$
and the family
$\mathcal{P}(\Z^n)$ of all lattice polytopes; namely, convex hulls of finitely many points of $\Z^n$
(see McMullen \cite{McM09}).
If $\mathcal{A}$ is a cancellative monoid (cancellative semigroup with identity element $0_\mathcal{A}$), then a function $Z:\mathcal{F}\to \mathcal{A}$
is called a {\it valuation} if the following holds: if $P\cup Q$ is convex for $P,Q\in\mathcal{F}$, then
\begin{equation}
\label{valuationdef}
Z(P\cup Q)+Z(P\cap Q)=Z(P)+Z(Q).
\end{equation}
We say that the valuation $Z:\mathcal{F}\to \mathcal{A}$ is {\it simple} if $Z(P)=0_\mathcal{A}$ for any $P\in\mathcal{F}$ with ${\rm dim}\,P\leq n-1$. Typically, one would consider valuations intertwining with some natural group actions, as we will shortly see.

While the idea of valuations on convex polytopes played a crucial role in Dehn's solution of Hilbert's Third Problem already around 1900, after sporadic results, the systematic study of valuations only started with Hadwiger's celebrated characterization of the intrinsic volumes as the basis of the space of continuous isometry invariant valutions from 1957. For the  breathtaking developments of the last seven decades, see for example the monograph Alesker \cite{Ale18}, and the survey papers
Alesker \cite{Ale17}, Ludwig \cite{Lud23} and Ludwig, Mussnig \cite{LuM23}. 
 The theory of valuations on lattice polytopes has been flourishing since the classical paper by Betke, Kneser \cite{BeK85} in 1985 characterizing unimodular invariant valuations on lattice polytopes, see for example, B\"or\"oczky, M. Ludwig \cite{BoL17,BoL19}, 
Jochemko, Sanyal \cite{JoS17,JoS18} and Ludwig, Silverstein \cite{LuS17}.

To state the result stimulating our research, let $L^1_c(\R^n)$ denote the family
of Lebesgue integrable real functions with compact support on $\R^n$. In addition, for measurable $\Omega\subset\R^n$, let $\mathcal{M}(\Omega)$ denote the family of real valued measurable functions on $\Omega$ where measurable means Lebesgue measurable without stating it in this paper. Actually, all the arguments and statements in this paper apply to the case if we fix a finite dimensional real vector space $\mathcal{V}$, and $\mathcal{M}(\Omega)$ denotes the family of measurable functions $\Omega\to \mathcal{V}$, including the case of complex valued measurable functions. However, for simplicity, we just discuss real valued functions.
For $f\in L^1_c(\R^n)$, its Laplace transform is
$$
\mathcal{L}f(x)=\int_{\R^n}e^{-\langle x,y\rangle}f(y)\,dy.
$$
Li, Ma \cite{LiM17} extended the definition of the Laplace transform to  a  convex compact set $K$ by applying $\mathcal{L}$
to the characteristic function; namely,
$$
\mathcal{L}K=\mathcal{L}\mathbf{1}_K=\int_Ke^{-\langle x,y\rangle}\,dy.
$$

Inspired by the properties of the Laplace transform
applied to the characteristic functions of compact convex sets,
Li, Ma \cite{LiM17} considered valuations  that are
 translatively exponential and ${\rm GL}(n,\R)$ covariant. For a subgroup $G\subset {\rm GL}(n,\R)$,
we say that a valuation $Z:\mathcal{K}^n\to \mathcal{M}(\R^n)$ is translatively exponential and $G$ covariant
if for any $x\in\R^n$, we have
\begin{eqnarray}
\label{unimod-inv-rn}
Z(\Phi K)(x)&=&|\det \Phi|\cdot Z( K)(\Phi^Tx) \mbox{ \ for any $\Phi\in G$;}\\
\label{log-trans-inv-rn}
Z( K+z)(x)&=&e^{-\langle z,x\rangle} Z(K)(x) \mbox{ \ for any $z\in\R^n$},
\end{eqnarray}
respectively.
The Hausdorff distance $\delta_H(C,K)$ between compact convex sets $K,C\subset\R^n$ is the minimal $r\geq 0$ such that
$K\subset C+r B^n$ and $C\subset K+r B^n$ where $B^n$ is the unit ball in $\R^n$ centered at the origin.
We say that a valuation $Z:\mathcal{K}^n\to \mathcal{M}(\R^n)$ is continuous if
for any convergent sequence $K_i\to K$ of compact convex sets with respect to the Hausdorff distance,
$Z(K_i)$ tends pointwise to $Z(K)$. In addition, $C(\R^n)$ denotes the space of real continuous functions on $\R^n$.

\begin{theo}[Li, Ma \cite{LiM17}]
Any continuous translatively exponential and ${\rm GL}(n,\R)$ covariant valuation $Z:\mathcal{K}^n\to C(\R^n)$
is of the form $Z=c\mathcal{L}$ for a constant $c\in\R$.
\end{theo}

Even if - as we will shortly see - there is an abundance of  translatively exponential and ${\rm GL}(2,\Z)$ covariant 
valuations $Z:\mathcal{P}(\Z^2)\to \mathcal{M}(\R^n)$ on lattice polytopes, we still believe that it is not so
for  translatively exponential and ${\rm GL}(n,\R)$ covariant valuations $Z:\mathcal{P}^n\to \mathcal{M}(\R^n)$ on polytopes in $\R^n$.

\begin{conj}
Any  translatively exponential and ${\rm GL}(n,\R)$ covariant valuation $Z:\mathcal{P}^n\to \mathcal{M}(\R^n)$
is of the form $Z=c\mathcal{L}$ for a constant $c\in\R$.
\end{conj}

In this paper, we consider some valuation $Z$ on lattice polygons (of $\Z^2$) with values in the space of measurable functions
that satisfy  \eqref{log-trans-inv}  and \eqref{unimod-inv}: If $x\in\R^2$, then
\begin{eqnarray}
\label{log-trans-inv}
Z( P+z)(x)&=&e^{\langle z,x\rangle} Z( P)(x) \mbox{ \ for any $z\in\Z^2$;}\\
\label{unimod-inv}
Z(\Phi P)(x)&=&Z( P)(\Phi^Tx) \mbox{ \ for any $\Phi\in {\rm GL}(2,\Z)$}
\end{eqnarray}
where \eqref{unimod-inv} is a restriction of \eqref{unimod-inv-rn} to ${\rm GL}(2,\Z)$ as $|\det \phi|=1$ for $\Phi\in {\rm GL}(2,\Z)$.
One example is the ``positive Laplace transform"
$$
\mathcal{L}_+(P)(x)=\int_{\R^n}e^{\langle x,y\rangle}\mathbf{1}_P(y)\,dy=\int_Pe^{\langle x,y\rangle}\,dy.
$$
We observe that $x\mapsto Z(P)(x)$ satisfies \eqref{log-trans-inv} if and only if
$x\mapsto Z(P)(-x)$ satisfies \eqref{log-trans-inv-rn}.
We note that Freyer, Ludwig, Rubey \cite{FMR} characterized the translatively exponential and ${\rm GL}(2,\Z)$ covariant valuations with values in the space of formal power series in two variables. 

The main goal of this paper is to characterize translatively exponential  and ${\rm GL}(2,\Z)$ covariant
 valuations $Z:\mathcal{P}(\Z^2)\to \mathcal{M}(\R^2)$.
We write $e_1,e_2$ to denote the orthonormal basis of $\R^2$ also generating $\Z^2$,
and set $T=[e_1,e_2,o]$ where $o=(0,0)$ stands for the origin, and $[x_1,\ldots,x_k]$ stands for the convex hull of
$x_1,\ldots,x_k\in \R^2$. In the formulas below,
$$
\frac{e^t-1}{t}\mbox{ is identified with }\sum_{n=0}^\infty\frac{t^n}{(n+1)!};
$$
namely, when we write $\varphi(t)=\frac{e^t-1}{t}$, we mean the positive analytic function $\varphi(t)=\sum_{n=0}^\infty\frac{t^n}{(n+1)!}$
on $\R$ satisfying $\varphi(0)=1$.
For the golden ratio $\tau=\frac{\sqrt{5}+1}2$, we set
$$
\widetilde{\Omega}_{2}=\{(\tau s,s):1\leq s<\tau\}\cup \{(x_1,x_2):0\leq x_1\leq x_2\}.
$$

\begin{figure}
    \centering
    \includegraphics[width=0.3\linewidth]{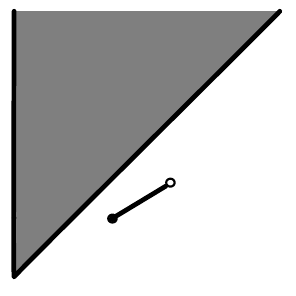}
    \caption{The domain $\widetilde{\Omega}_2$.}
    \label{fig:omega2}
\end{figure}

First, we characterize simple  translatively exponential  and ${\rm GL}(2,\Z)$ covariant
 valuations where a valuation $Z:\mathcal{P}(\Z^2)\to \mathcal{M}(\R^2)$ is {\it simple} if
$Z(P)=0$ for any  $P\in\mathcal{P}(\Z^2)$ with ${\rm dim}\,P\leq 1$.

\begin{theo}
\label{simpleGLcovariant}
Simple translatively exponential  and ${\rm GL}(2,\Z)$ covariant
 valuations $Z:\mathcal{P}(\Z^2)\to \mathcal{M}(\R^2)$ are parametrized  uniquely by functions in
$\mathcal{M}(\widetilde{\Omega}_{2})$ in the following way.
\begin{description}
\item{(i)}
For any measurable function $\tilde{\varrho}:\widetilde{\Omega}_{2}\to\R$, there exists
a unique measurable extension $\varrho:\,[0,\infty)^2\to\R$ satisfying $\tilde{\varrho}=\varrho|_{\widetilde{\Omega}_{2}}$ and
\begin{equation}
\label{rhoformula}
(2x+y) \varrho(x,y)=(x+y)\varrho(x,x+y)+x \varrho(x+y,x)
\end{equation}
for $x,y\geq 0$, and a unique  translatively exponential  and ${\rm GL}(2,\Z)$ covariant
simple valuation $Z:\mathcal{P}(\Z^2)\to \mathcal{M}(\R^2)$  such that
$Z(T)(0,0)=\tilde{\varrho}(0,0)$, and
 $f_2=Z(T)$ satisfies
\begin{equation}
\label{simpleZrho}
f_2(x,y)=\frac{e^x}{y}\cdot \frac{e^{y-x}-1}{y-x}\cdot \varrho(y-x,x)-
\frac{1}y\cdot \frac{e^{x}-1}{x}\cdot \varrho(x,y-x)
\end{equation}
for $0\leq x\leq y$ and $y>0$.
\item{(ii)} For any simple translatively exponential  and ${\rm GL}(2,\Z)$ covariant
 valuation $Z:\mathcal{P}(\Z^2)\to \mathcal{M}(\R^2)$, there exists some
measurable $\varrho:\,[0,\infty)^2\to\R$ satisfying \eqref{rhoformula}
and \eqref{simpleZrho}.
\end{description}
\end{theo}
\noindent{\bf Remark.} The ``positive Laplace transform" $\mathcal{L}_+$  corresponds to the constant one function $\varrho\equiv 1$. On the other hand,  Freyer, Ludwig, Rubey \cite{FMR} determined all translatively exponential  and ${\rm GL}(2,\Z)$ covariant valuations with values in the space of formal power series in two variables. 
It follows that there exists an abundance of corresponding simple valuations with values in the space of 
 analytic functions
that are different from the positive Laplace transform, and we exhibit one in Example~\ref{analytic-example}.\\

We write $|X|$ to denote the Lebesgue measure of a measurable set $X\subset \R^2$.
We need the following consequence of the fact that the linear action of ${\rm SL}(2,\Z)$ on $\R^2$ is
ergodic (see Section~\ref{secSL2Zergodic}).
Proposition~\ref{fSL2Zinvariant} about Lebesgue measurable ${\rm SL}(2,\Z)$ invariant functions is used in the statement of Theorem~\ref{ExpValGL2Z} characterizing translatively exponential  and ${\rm GL}(2,\Z)$ covariant
valuations.

\begin{prop}
\label{fSL2Zinvariant}
If $f\in \mathcal{M}(\R^2)$ is invariant under ${\rm SL}(2,\Z)$, then there exists  a constant $c\in\R$
such that $f(x)=c$ for almost everywhere $x\in\R^2$.
\end{prop}
\noindent{\bf Remark. } In particular, any function $f\in \mathcal{M}(\R^2)$ invariant under ${\rm GL}(2,\Z)$ can be constructed in the following way from a $c\in\R$ and an $X\subset \R^2$ with $|X|=0$. Writing $X_0$ to denote the image of $X$ under the action
of ${\rm GL}(2,\Z)$, we have $|X_0|=0$,
and associating
an arbitrary $c_\mathcal{O}\in\R$ to any orbit $\mathcal{O}$ of ${\rm GL}(2,\Z)$ intersecting $X$, we
define
$f(x)=c$ for $x\in\R^2\backslash X_0$,  and  $f(x)=c_\mathcal{O}$ if $x\in\mathcal{O}$ for
an orbit $\mathcal{O}$ of ${\rm GL}(2,\Z)$ intersecting $X$.\\

In Theorem~\ref{ExpValGL2Z}, we use
$$
\widetilde{\Omega}_1=\{(0,x_2):x_2\geq 0\}\cup
\left\{(x_1,x_2):0\leq x_2\leq \mbox{$\frac12\,$}x_1\right\}.
$$
\begin{figure}
    \centering
    \includegraphics[width=0.3\linewidth]{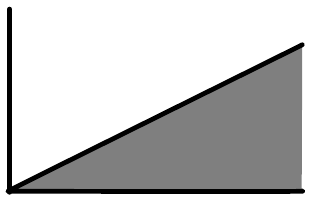}
    \caption{The domain $\widetilde{\Omega}_1$.}
    \label{fig:omega1}
\end{figure}

The space of ${\rm GL}(2,\Z)$ invariant measurable functions in $\R^2$ - that are characterized in Proposition~\ref{fSL2Zinvariant} -, is denoted by $\mathcal{M}(\R^2)^{{\rm GL}(2,\Z)}$.

\begin{theo}
\label{ExpValGL2Z}
Translatively exponential  and ${\rm GL}(2,\Z)$ covariant
valuations $Z:\mathcal{P}(\Z^2)\to \mathcal{M}(\R^2)$ are parameterized uniquely
by $\mathcal{M}(\R^2)^{{\rm GL}(2,\Z)}$, $\mathcal{M}(\widetilde{\Omega}_1)$ and $\mathcal{M}(\widetilde{\Omega}_2)$ as follows:
\begin{description}
\item{(i)}
For any ${\rm GL}(2,\Z)$ invariant measurable function $f_0:\R^2\to\R$ (cf. Proposition~\ref{fSL2Zinvariant}),
any measurable functions $\tilde{f}_1:\widetilde{\Omega}_{1}\to\R$,
and $\tilde{\rho}:\widetilde{\Omega}_{2}\to\R$,  there exists a unique
translatively exponential and ${\rm GL}(2,\Z)$ covariant
valuation $Z:\mathcal{P}(\Z^2)\to \mathcal{M}(\R^2)$ such that $Z(\{o\})=f_0$,
$f_1|_{\widetilde{\Omega}_1}=\tilde{f}_1$
for $f_1=Z([o,e_1])$, and
\begin{equation}
\label{ZZ2Z1}
Z(T)(x,y)=f_2(x,y)+\mbox{$\frac12$}\,f_1(x,y)+
\mbox{$\frac12$}\,f_1(y,-x)
+\mbox{$\frac{e^{x}}2$}\,f_1(-x+y,-x)
\end{equation}
holds for the $f_2$ defined by \eqref{simpleZrho} where
$f_2=Z_2(T)$ for a simple 
translatively exponential ${\rm GL}(2,\Z)$ covariant
valuation $Z_2:\mathcal{P}(\Z^2)\to \mathcal{M}(\R^2)$
and
$\varrho$ is constructed 
from $\tilde \varrho$ as in Theorem~\ref{simpleGLcovariant} (i).

\item{(ii)} For a  translatively exponential and ${\rm GL}(2,\Z)$ covariant
valuation $Z:\mathcal{P}(\Z^2)\to \mathcal{M}(\R^2)$, $f_0=Z(\{o\})$ is a
measurable ${\rm GL}(2,\Z)$ invariant function and
\begin{equation}
\label{Zf_0p}
Z(\{p\})(z)=e^{\langle p,z\rangle}f_0(z)\mbox{ \ for $p\in\Z^2$ and $x\in\R^2$,}
\end{equation}
where there exists constant $c\in\R$ such that $f_0=c$ almost everywhere.

In addition, there exists
a simple translatively exponential  and ${\rm GL}(2,\Z)$ covariant
valuation $Z_2:\mathcal{P}(\Z^2)\to \mathcal{M}(\R^2)$ such that
$f_1=Z([o,e_1])$ and $f_2=Z_2(T)$ satisfy \eqref{ZZ2Z1}.

\end{description}
\end{theo}
\noindent{\bf Remark. } In particular, any  translatively exponential  and ${\rm GL}(2,\Z)$ covariant
valuation $Z:\mathcal{P}(\Z^2)\to \mathcal{M}(\R^2)$ can be represented as
$Z=Z_1+Z_2$ where $Z_2$ is a simple valuation and $Z_1$
is the valuation constructed in
Proposition~\ref{lowdim} (see Section~\ref{secExamples}) 
satisfying $Z_1(\{o\})=Z(\{o\})$ and  $Z_1([o,e_1])=Z([o,e_1])$.\\

Concerning the structure of the paper, properties of related group actions are discussed in Section~\ref{secGroupAction}, and
Section~\ref{secExamples} presents some fundamental examples. For a translatively exponential  and ${\rm GL}(2,\Z)$ covariant
valuation $Z:\mathcal{P}(\Z^2)\to \mathcal{M}(\R^2)$,
the algebraic properties characterizing
$Z(\{o\})=f_0$,
 $Z([o,e_1])=f_1$ and $Z_2(T)=f_2$ for the corresponding translatively exponential  and ${\rm GL}(2,\Z)$ covariant
simple valuation $Z_2:\mathcal{P}(\Z^2)\to \mathcal{M}(\R^2)$
are described in
Lemma~\ref{f0f1f2properties} and in Section~\ref{secsimplef2}.

The restrictions of  translatively exponential  and ${\rm GL}(2,\Z)$ covariant
valuations to points and lattice segments are characterized in Section~\ref{secSL2Zergodic} and Section~\ref{secf1}, and  Theorem~\ref{ExpValGL2Z} is actually proved in Section~\ref{secf1}.
The study of geometric properties of simple translatively exponential  and ${\rm GL}(2,\Z)$ covariant
valuations starts in Section~\ref{secf2step1}, and finally Section~\ref{secf2full} verifies Theorem~\ref{simpleGLcovariant} using the Fibonacci sequence.

\section{Action of affine lattice automorphisms}
\label{secGroupAction}

We write $\mathcal{G}(\Z^2)$ to denote the group of all affine transforms $\Xi x=\Phi x+a$ for
$\Phi\in {\rm GL}(2,\Z)$ and $a\in\Z^2$, and hence $\mathcal{G}(\Z^2)$
is the group of all affine automorphisms of $\Z^2$.

Let us define a $\mathcal{G}(\Z^2)$ action on  $\mathcal{M}(\R^2)$.
For $f\in \mathcal{M}(\R^2)$ and an affine $\Xi x=\Phi x+a$ with
$\Phi\in {\rm GL}(2,\Z)$ and $a\in\Z^2$,
we consider $\Xi\cdot f\in \mathcal{M}(\R^2)$ where
\begin{equation}
\label{AffineActionOnMeasurable}
(\Xi\cdot f)(x)=e^{\langle a,x\rangle} f(\Phi^Tx).
\end{equation}
Formula \eqref{AffineActionOnMeasurable} defines an action of $\mathcal{G}(\Z^2)$
because the identity map leaves any $f\in \mathcal{M}(\R^2)$ invariant, and for any $\aleph\in \mathcal{G}(\Z^2)$,
we have
\begin{equation}
\label{AffineActionOnMeasurableComposition}
(\aleph\circ \Xi)\cdot f=\aleph\cdot (\Xi\cdot f).
\end{equation}
To prove \eqref{AffineActionOnMeasurableComposition}, let
$\aleph x=\Psi x+b$ for
$\Psi\in {\rm GL}(2,\Z)$ and $b\in\Z^2$. Hence
$(\aleph\circ \Xi)(x)=\Psi\Phi x+\Psi a+b$ and
\begin{align*}
\Big((\aleph\circ \Xi)\cdot f\Big)(x)&=e^{\langle \Psi a+b,x\rangle} f(\Phi^T\Psi^Tx)=
e^{\langle b,x\rangle}e^{\langle  a,\Psi^T x\rangle}f(\Phi^T\Psi^Tx) \\
&=e^{\langle b,x\rangle}(\Xi\cdot f)(\Psi^T x)=\Big(\aleph\cdot (\Xi\cdot f)\Big)(x).
\end{align*}
We observe that this action of a $\Xi\in\mathcal{G}(\Z^2)$ satisfies that for a valuation $Z:\mathcal{P}(\Z^2)\to \mathcal{M}(\R^2)$, we have 
$\Xi\cdot Z(P)=Z(\Xi P)$ for $P\in \mathcal{P}(\Z^2)$.

The formula \eqref{AffineActionOnMeasurableComposition} readily yields the following property.

\begin{lemma}
\label{well-defined}
Let $P\in\mathcal{P}(\Z^2)$, and let
 $H\subset \mathcal{G}(\Z^2)$ be the subgroup of all $\Xi\in \mathcal{G}(\Z^2)$ such that $\Xi P=P$.

For a set $\{\Xi_i\}_{i\in I}$ of group generators of $H$, if $\xi(P)\in \mathcal{M}(\R^2)$ satisfies
$\Xi_i\cdot \xi(P)=\xi(P)$ for $i\in I$, then
setting $\xi(\aleph P)\in \mathcal{M}(\R^2)$ via $\xi(\aleph P)=\aleph\cdot \xi(P)$
for any $\aleph\in \mathcal{G}(\Z^2)$ is well-defined, and satisfies
$$
\xi\Big(\Gamma(\aleph P)\Big)=\Gamma\cdot \xi(\aleph P) \mbox{ \ for any $\Gamma\in \mathcal{G}(\Z^2)$}.
$$
\end{lemma}

For a $P\in\mathcal{P}(\Z^2)$, Lemma~\ref{StabilizerGenertors} is used to determine the subgroup
$\mathcal{G}(\Z^2)$ fixing $P$. For a $\Xi\in \mathcal{G}(\Z^2)$
with $\Xi x=\Phi x+a$, $\Phi\in {\rm GL}(2,\Z)$ and $a\in \Z^2$, we set $\det \Xi=\det \Phi$. Since
$\det: \mathcal{G}(\Z^2)\to\{-1,1\}$ is a group homomorphism, and a vertex of $P$ is mapped to a vertex of the image by
an affine transform (where endpoints of segments are understood as vertices), we deduce the following. We call
a $\Xi\in \mathcal{G}(\Z^2)$ orientation preserving if $\det\Xi=1$.

\begin{lemma}
\label{StabilizerGenertors}
Let $P\in\mathcal{P}(\Z^2)$ such that $o$ is a vertex, let $H\subset \mathcal{G}(\Z^2)$
 be the subgroup of all $\Xi\in  \mathcal{G}(\Z^2)$
 such that $\Xi P=P$, and let $H_+\subset H$ be the subgroup of orientation preserving transformations in $H$ .
\begin{description}
\item{(i)} $\Xi(v)$ is a vertex of $P$ for any vertex $v$ of $P$ and  $\Xi\in H$;
\item{(ii)} if there exists $\widetilde{\Xi}\in  \mathcal{G}(\Z^2)$ with $\widetilde{\Xi} P=P$ and  $\det \widetilde{\Xi}=-1$,
then $H_+$ and $\widetilde{\Xi}$ generate $H$;
\item{(iii)} if $P$ is a polygon and $\Xi(o)=\Xi'(o)$ for
$\Xi,\Xi'\in H_+$, then $\Xi=\Xi'$.
\end{description}
\end{lemma}

For a subgroup $G\subset \mathcal{G}(\Z^2)$, we say that
 a valuation   $Z:\mathcal{P}(\Z^2)\to \mathcal{M}(\R^2)$ is $G$
covariant if
\begin{equation}
\label{Gcontravariant}
Z(\Xi P)(x)=\Big(\Xi\cdot Z(P)\Big)(x)\mbox{ \ for any
$\Xi\in G$, $P\in \mathcal{P}(\Z^2)$ and $x\in\R^2$.}
\end{equation}
In particular, $Z$ is $\mathcal{G}(\Z^2)$ covariant if and only if
$Z$ is translatively exponential  ({\it cf.} \eqref{log-trans-inv})
and ${\rm GL}(2,\Z)$ covariant ({\it cf.} \eqref{unimod-inv}).

Next we show that a
$\mathcal{G}(\Z^2)$ covariant
valuation $Z:\mathcal{P}(\Z^2)\to \mathcal{M}(\R^2)$
 is determined by its value at $\{o\}$, $[o,e_1]$ and $T$.
We note that McMullen \cite{McM09} proved the inclusion-exlusion principle for any valuation
$Z:\mathcal{P}(\Z^m)\to \mathcal{A}$ where $\mathcal{A}$ is a cancellative semigroup; namely, if
$P$ is a $d$-dimensional lattice polytope, $1\leq d\leq m$, and $P=Q_1\cup\ldots\cup Q_k$
for $d$-dimensional lattice polytopes $Q_1,\ldots, Q_k$, $k\geq 3$, such that the non-empty intersection
of any subfamily of $Q_1,\ldots, Q_k$ is a lattice polytope, then
\begin{equation}
\label{IncExc}
Z(P)=\sum_{i=1}^k\sum_{1\leq j_1<\ldots<j_i\leq k}(-1)^{i-1}Z(Q_{j_1}\cap\ldots Q_{j_i})
\end{equation}
where we set $Z(\emptyset)=0$.

We say that a lattice segment $[p,q]$ for $p,q\in \Z^2$ is a {\it primitive lattice segment} if the only lattice points it contains are its endpoints,
and a lattice triangle is an {\it empty triangle} if the only lattice points it contains are its vertices.

\begin{prop}
\label{ExpValSL2Z}
For  any  $\mathcal{G}(\Z^2)$ covariant
valuations $Z,Z':\mathcal{P}(\Z^2)\to \mathcal{M}(\R^2)$,
\begin{description}
\item{(i)} if $Z(\{o\})=Z'(\{o\})$ and $Z([o,e_1])=Z'([o,e_1])$, then
$Z-Z'$ is a simple $\mathcal{G}(\Z^2)$ covariant
valuation;
\item{(ii)}
if $Z(\{o\})=Z'(\{o\})$, $Z([o,e_1])=Z'([o,e_1])$ and $Z(T)=Z'(T)$,
then $Z=Z'$.
\end{description}
\end{prop}
\proof Let $Z(\{o\})=Z'(\{o\})$ and $Z([o,e_1])=Z'([o,e_1])$. For any $p\in\Z^2$, \eqref{log-trans-inv}
yields that
$$
Z(\{p\})(x)=e^{\langle x,p\rangle}Z(\{o\})(x)=e^{\langle x,p\rangle}Z'(\{o\})(x)=Z'(\{p\})(x).
$$

For a primitive lattice segment $[p,q]$, $p\neq q\in \Z^2$, there exists
$\Phi\in{\rm SL}(2,\Z)$ such that $[p,q]=p+\Phi[o,e_1]$; therefore,
$Z([p,q])=Z'([p,q])$ by
\eqref{log-trans-inv} and \eqref{unimod-inv}. Now a lattice segment $s$ containing $k\geq 3$ lattice points can be subdivided into $k-1$ primitive lattice segments, and hence the Inclusion-Exclusion principle \eqref{IncExc}, together with the just established property that $Z$ and $Z'$  agree on lattice points and primitive segments, yields that $Z(s)=Z'(s)$.
In turn, we deduce (i).

Let us assume that in addition, $Z(T)=Z'(T)$.
For an empty lattice triangle
$T=[v_1,v_2,v_3]$, the vectors $v_2-v_1$ and $v_3-v_1$ form a basis of $\Z^2$; therefore,
there exists a $\Phi\in{\rm SL}(2,\Z)$ such that $T=v_1+\Phi T$, which in turn yields that
$Z(T)=Z'(T)$ by
\eqref{log-trans-inv} and \eqref{unimod-inv}. Finally, a two-dimensional lattice polygon $P$ can be written
as $P=T_1\cup\ldots\cup T_k$ for empty triangles
 $T_1,\ldots,T_k$ such that $T_i\cap T_j$ is either empty, or a common vertex, or a common side. Thus
the Inclusion-Exclusion principle \eqref{IncExc} together with the property that $Z$ and $Z'$ agree on lattice points, lattice segments, and empty triangles yields that $Z(P)=Z'(P)$.
\endproof

It follows from Proposition~\ref{ExpValSL2Z} that
in order to characterize a translatively exponential  and
 ${\rm GL}(2,\Z)$ covariant (or equivalently, $\mathcal{G}(\Z^2)$ covariant) valuation $Z$ on lattice polygons,
all we need to characterize are
\begin{equation}
\label{f0f1f2}
Z(\{o\})=f_0\mbox{ and }
Z([o,e_1])=f_1\mbox{ and }
Z(T)=f_2.
\end{equation}

Let us establish some algebraic properties of $f_0$ and $f_1$
in \eqref{f0f1f2}.
According to Trott \cite{Tro62},
 ${\rm GL}(2,\Z)$ as a group is generated by
$\left\{
\left[
\begin{array}{cc}
0&1\\
1&0
\end{array}\right],\mbox{ }
\left[
\begin{array}{cc}
1&1\\
0&1
\end{array}\right]
\right\}
$, but we do not use this result directly.

\begin{lemma}
\label{f0f1f2properties}
Let
$Z:\mathcal{P}(\Z^2)\to\mathcal{M}(\R^2)$ be a $\mathcal{G}(\Z^2)$ covariant
valuation. Then
$Z(\{o\})=f_0$ and $Z([o,e_1])=f_1$ satisfy the following properties:
\begin{eqnarray}
\label{f0GL2ZSL2Z0}
f_0(\Phi x)&=&f_0(x) \mbox{ \ for any $\Phi\in {\rm GL}(2,\Z)$ and $x\in \R^2$};\\
\label{f1shift}
f_1(-x_1,-x_2)&=& e^{-x_1}f_1(x_1,x_2)\mbox{ \ for $(x_1,x_2)\in\R^2$};\\
\label{f1111}
f_1(x_1,x_2)&=& f_1(x_1,x_1+x_2)\mbox{ \ for $(x_1,x_2)\in\R^2$};\\
\label{f1extra}
f_1(x_1,x_2)&=& f_1(x_1,-x_2)\mbox{ \ for $(x_1,x_2)\in\R^2$}.
\end{eqnarray}
\end{lemma}
\proof The formula \eqref{f0GL2ZSL2Z0} follows from  \eqref{unimod-inv}.
The property \eqref{f1shift} follows from \eqref{log-trans-inv} and the
relation
$\left[
\begin{array}{cc}
-1&0\\
0&-1
\end{array}\right][o,e_1]=[o,e_1]-e_1$. Next \eqref{f1111} says that $Z([o,e_1])$ is covariant under
$\left[
\begin{array}{cc}
1&1\\
0&1
\end{array}\right]$, and \eqref{f1extra} is the consequence of the fact that
 $Z([o,e_1])$ also is covariant under
$\left[
\begin{array}{cc}
1&0\\
0&-1
\end{array}\right]$.
\endproof

\section {Some examples of $\mathcal{G}(\Z^2)$ covariant valuations on lattice polygons}
\label{secExamples}

In this section, we construct two fundamental examples of $\mathcal{G}(\Z^2)$ covariant valuations on $\mathcal{P}(\Z^2)$. The first is the "positive Laplace transform" that is the basic example of simple valuations satisfying  \eqref{log-trans-inv} and \eqref{unimod-inv}, and the second in Propositions~\ref{lowdim} extends a "would be" valuation defined on lattice points and lattice segments.

\begin{example}[Positive Laplace transform]
For any lattice polygon $P$, we define
\begin{equation}
\label{positiveLaplace}
\mathcal{L}_+(P)(x)=\int_Pe^{\langle x,y\rangle}\,dy,
\end{equation}
thus $\mathcal{L}_+$ is a $\mathcal{G}(\Z^2)$ covariant simple valuation.

We note that \eqref{positiveLaplace} yields that
\begin{equation}
\label{LaplaceTformulaxt}
\mathcal{L}_+(T)(x_1,x_2)=\frac{x_1e^{x_2}-x_2 e^{x_1}+x_2-x_1}{x_1x_2(x_2-x_1)};
\end{equation}
or in other words, we have
\begin{equation}
\label{LaplaceTformula}
\mathcal{L}_+(T)(x_1,x_1+x_2)=\frac1{x_1+x_2}\left(e^{x_1}\cdot \frac{e^{x_2}-1}{x_2}- \frac{e^{x_1}-1}{x_1}\right).
\end{equation}
\end{example}

\begin{prop}
\label{lowdim}
For $f_0,f_1\in\mathcal{M}(\R^2)$ satisfying  \eqref{f0GL2ZSL2Z0}, \eqref{f1shift}, \eqref{f1111} and
 \eqref{f1extra}, there exists
a  $\mathcal{G}(\Z^2)$ covariant valuation $Z_1:\mathcal{P}(\Z^2)\to\mathcal{M}(\R^2)$ satisfying
$Z_1(\{o\})=f_0$ and $Z_1([o,e_1])=f_1$ and
$$
Z_1(T)(x_1,x_2)=\mbox{$\frac12$}\,f_1(x_1,x_2)+\mbox{$\frac12$}\,f_1(x_2,-x_1)
+\mbox{$\frac{e^{x_1}}2$}\,f_1(-x_1+x_2,-x_1).
$$
\end{prop}
\proof
For points, we define
\begin{equation}
\label{Z0exa}
Z_1(\{z\})(x)=e^{\langle x,z\rangle}f_0(x) \mbox{ \ for $z\in\Z^2$ and $x\in \R^2$}.
\end{equation}
As ${\rm GL}(2,\Z)\subset \mathcal{G}(\Z^2)$ is the stabilizer subgroup of $\{o\}$, \eqref{f0GL2ZSL2Z0} and Lemma~\ref{well-defined} yield that
$Z_1$ as defined in \eqref{Z0exa} is well-defined and is $\mathcal{G}(\Z^2)$ covariant on lattice points.

Turning to segments, we write $H$ to denote the stabilizer subgroup of $[o,e_1]$ in 
$\mathcal{G}(\Z^2)$, $H_+\subset H$ to denote the orientation preserving stabilizers,  and
$H_0\subset {\rm SL}(2,\Z)$ to denote the subgroup of all $\Phi\in{\rm SL}(2,\Z)$ with $\Phi e_1=e_1$.
Now $H_0$ consists of the matrices of the form $\left[
\begin{array}{cc}
1&a\\
0&1
\end{array}\right]$ for an $a\in\Z$, and hence it it is generated by
$\Phi_0=\left[
\begin{array}{cc}
1&1\\
0&1
\end{array}\right]$. Since
$\Xi_1[o,e_1]=[o,e_1]$ for $\Xi_1x=-x+e_1\in H_+\backslash H_0$, and
$H_0\subset H_+$ is a subgroup of index $2$ by Lemma~\ref{StabilizerGenertors} (i),
we deduce that $H_+$ is generated by $\Phi_0$ and $\Xi_1$. It follows from
Lemma~\ref{StabilizerGenertors} (ii) that $H$ is generated by $\Phi_0$, $\Xi_1$
and $\Phi_{-}=\left[
\begin{array}{cc}
1&0\\
0&-1
\end{array}\right]$.

For any primitive lattice segment $[p,q]$ there exists $\Xi\in \mathcal{G}(\Z^2)$ such that
$\Xi[o,e_1]=[p,q]$; namely, take $\Xi x=\Phi x+p$ where
$q-p=\Phi e_1$ for $\Phi\in {\rm SL}(2,\Z)$. Since
$H$ is generated by $\Phi_0$, $\Xi_1$
and $\Phi_{-}$ where the condition \eqref{f1extra} on $f_1$ corresponds to the $\Phi_{-}$
invariance of $[o,e_1]$, it follows from Lemma~\ref{well-defined} that
if $[p,q]$ is a primitive lattice segment, then
  we may define
\begin{equation}
\label{Z1primitivedef}
Z_1([p,q])=\Xi\cdot f_1
\end{equation}
for any $\Xi\in \mathcal{G}(\Z^2)$ with $\Xi[o,e_1]=[p,q]$.
  Lemma~\ref{well-defined} also yields that
$Z_1$ is $\mathcal{G}(\Z^2)$ covariant on lattice segments.

Now a general lattice segment $[a,b]$ for $a\neq b$, $a,b\in\Z^2$ is divided into primitive lattice segments
by the lattice points contained in $[a,b]$ different from $a,b$, and hence,
using the abbreviation "pls"  for primitive lattice segments, we define
\begin{equation}
\label{Z1segmentdef}
Z_1([a,b])=\sum_{\substack{s\subset[a,b]\\ s\;{\rm pls}}}Z_1(s)-\sum_{\substack{z\in [a,b]\cap \Z^2\\ z\neq a,b}}Z_1(\{z\}).
\end{equation}
Since $Z_1$ is translatively exponential and ${\rm GL}(\Z^2)$ covariant on points and primitive lattice segments, $Z_1$ is $\mathcal{G}(\Z^2)$ covariant on lattice segments and lattice points.

In addition, if $P$ is a lattice polygon, then we define
\begin{equation}
\label{Z1polygondef}
Z_1(P)=\sum_{\substack{s\subset\partial P\\ s\;{\rm pls}}}\mbox{$\frac12\,$}Z_1(s)+\sum_{z\in ({\rm int} P)\cap \Z^2}Z_1(\{z\})
\end{equation}
where the second sum is zero if $({\rm int} P)\cap \Z^2=\emptyset$.
 $Z_1$ is readily $\mathcal{G}(\Z^2)$ covariant also on polygons.
In addition, the (primitive) edges of $T$ are $[o,e_1]$, $\Phi_1[o,e_1]$ and $e_1+\Phi_2[o,e_1]$
for
$\Phi_1=\left[
\begin{array}{cc}
0&-1\\
1&0
\end{array}\right]\in{\rm SL}(2,\Z)$ and
$\Phi_2=\left[
\begin{array}{cc}
-1&-1\\
1&0
\end{array}\right]\in{\rm SL}(2,\Z)$, 
thus \eqref{Z1primitivedef} and \eqref{Z1polygondef} yield
$$
Z_1(T)(x,y)=\mbox{$\frac12$}\,f_1(x,y)+\mbox{$\frac12$}\,f_1(y,-x)
+\mbox{$\frac{e^{x}}2$}\,f_1(-x+y,-x).
$$

Therefore, all we have to prove is that the function $Z_1$ on $\mathcal{P}(\Z^2)$ defined by
\eqref{Z0exa}, \eqref{Z1primitivedef}, \eqref{Z1segmentdef} and \eqref{Z1polygondef} is a valuation; namely, if
$P\cup Q$ is convex for $P,Q\in\mathcal{P}(\Z^2)$ , then
\begin{equation}
\label{Z1valuation}
Z_1(P)+Z_1(Q)=Z_1(P\cap Q)+Z_1(P\cup Q).
\end{equation}
If $P\subset Q$ or $Q\subset P$, then  \eqref{Z1valuation} readily holds. Therefore, besides that
$P\cup Q$ are convex, we assume that either
$P$ and $Q$ are lattice segments not containing each other, or
$P$ and $Q$ are lattice polygons not containing each other.\\

\noindent{\bf Case 1} {\it ${\rm dim}(P\cap Q)<{\rm dim}\,P={\rm dim}\,Q$ }\\
In other words, either $P$ and $Q$  are collinear lattice segments whose intersection is a common endpoint, or 
$P$ and $Q$  are  lattice polygons whose intersection is a common side, and in addition, $P\cup Q$ is convex.
In the first case \eqref{Z1segmentdef}, and in the second case \eqref{Z1segmentdef} and \eqref{Z1polygondef} directly yield \eqref{Z1valuation}.\\

\noindent{\bf Case 2} {\it ${\rm dim}(P\cap Q)={\rm dim}\,P={\rm dim}\,Q\geq 1$ }\\
Let us introduce a notion of multiplicity of a primitive lattice segment $s$ and a $z\in\Z^2$ with respect to a
general lattice segment $[a,b]$ for $a\neq b$, $a,b\in\Z^2$ and a lattice polygon $N$ where we set $(a,b)=[a,b]\backslash\{a,b\}$:
\begin{align*}
m_{[a,b]}(s)&=1\mbox{ if }s\subset[a,b],\mbox{ \ \ and }
m_{[a,b]}(s)=0\mbox{ otherwise; }\\
m_{[a,b]}(z)&=1\mbox{ if }z\in(a,b),\mbox{ \ \ and }
m_{[a,b]}(z)=0\mbox{ otherwise; }\\
m_{N}(s)&=1\mbox{ if }s\subset\partial N,\mbox{ \ \ and }
m_{N}(s)=0\mbox{ otherwise; }\\
m_{N}(z)&=1\mbox{ if }z\in{\rm int}\,N,\mbox{ \ \ and }
m_{N}(z)=0\mbox{ otherwise. }
\end{align*}
In particular, 
\begin{align}
\label{Z1segmentdef0}
Z_1([a,b])&=\sum_{s\;{\rm pls}}m_{[a,b]}(s)\cdot Z_1(s)-\sum_{z\in\Z^2}m_{[a,b]}(z)\cdot Z_1(\{z\})  \\
\label{Z1polygondef0}
Z_1(N)&=\sum_{s\;{\rm pls}}m_{N}(s)\cdot \mbox{$\frac12\,$}Z_1(s)+\sum_{z\in\Z^2}m_{N}(z)\cdot Z_1(\{z\}).
\end{align}
If $P,Q\in\mathcal{P}(\Z^2)$ satisfy that $P\cup Q$ is convex  and ${\rm dim}(P\cap Q)={\rm dim}\,P={\rm dim}\,Q\geq 1$, and if $s$ is a primitive lattice segment  and $z\in\Z^2$, then we claim that
\begin{align}
\label{multiplicity-additive-s}
m_{P}(s)+m_{Q}(s)&=
m_{P\cap Q}(s)+m_{P\cup Q}(s)\\
\label{multiplicity-additive-z}
m_{P}(z)+m_{Q}(z)&=
m_{P\cap Q}(z)+m_{P\cup Q}(z).
\end{align}
If ${\rm dim}(P\cap Q)={\rm dim}\,P={\rm dim}\,Q= 1$, then
$({\rm relint}\,P)\cup ({\rm relint}\,Q)={\rm relint}(P\cup Q)$ (where ${\rm relint}\,P$ is the open segment determined by $P$), and hence
$z\in ({\rm relint}\,P)\cap ({\rm relint}\,Q)$ if and only if 
$z\in {\rm relint}\,P$ and $z\in {\rm relint}\,Q$. Since 
$s\subset P\cap Q$ for a the primitive lattice segment $s$ if and only if $s\subset P$ and $s\subset Q$, we conclude \eqref{multiplicity-additive-s} and \eqref{multiplicity-additive-z}.

Let ${\rm dim}(P\cap Q)={\rm dim}\,P={\rm dim}\,Q= 2$, and hence
\begin{equation}
\label{no-separating-line}
\mbox{no line separates $P$ and $Q$.}
\end{equation}
For \eqref{multiplicity-additive-z}, if $x\in {\rm int}(P\cup Q)$ and $x\in \partial P$, then there exists a supporting line $\ell$ to $P$ at $x$. Now any small semicircle centered at $x$ and separated from $P$ by $\ell$ is contained in $Q$, and 
\eqref{no-separating-line} yields that $x\in {\rm int}\,Q$. The similar argument if $x\in {\rm int}(P\cup Q)$ and $x\in \partial Q$ shows that
$({\rm int}\,P)\cup ({\rm int}\,Q)={\rm int}(P\cup Q)$,
which in turn implies \eqref{multiplicity-additive-z}.
For \eqref{multiplicity-additive-s}, we deduce from the convexity of $P\cup Q$ that each vertex of $P\cap Q$ is a vertex of $P$ or $Q$ (see for example McMullen \cite{McM09}), and hence the primitive lattice segment $s$ is contained in $\partial(P\cup Q)$ or $\partial(P\cap Q)$ if and only if $s$ is contained in $\partial\,P$ or $\partial \,Q$.
We deduce from \eqref{no-separating-line} that if $s\subset(\partial\,P)\cap (\partial \,Q)$, then $s\subset \partial(P\cup Q)$ and $s\subset\partial(P\cap Q)$; therefore, 
\eqref{multiplicity-additive-s} holds, as well.

Combining 
\eqref{Z1segmentdef0},
\eqref{Z1polygondef0},
\eqref{multiplicity-additive-s} and \eqref{multiplicity-additive-z}
yields \eqref{Z1valuation}, and in turn Proposition~\ref{lowdim}.
\endproof

Proposition~\ref{ExpValSL2Z} (i) and Proposition~\ref{lowdim} yield the following.

\begin{coro}
\label{ReductionToSimple}
For any  $\mathcal{G}(\Z^2)$ covariant
valuation $Z:\mathcal{P}(\Z^2)\to \mathcal{M}(\R^2)$, let
$Z_1$ be the $\mathcal{G}(\Z^2)$ covariant
valuation satisfying $Z_1(\{o\})=Z(\{o\})$ and
$Z_1([o,e_1])=Z([o,e_1])$
constructed in Proposition~\ref{lowdim}. Then $Z_2=Z-Z_1$
is a $\mathcal{G}(\Z^2)$ covariant simple
valuation.
\end{coro}

\section{Algebraic properties of simple valuations}
\label{secsimplef2}

Having Lemma~\ref{f0f1f2properties}
and Corollary~\ref{ReductionToSimple} at hand, an algebraic characterization
of simple $\mathcal{G}(\Z^2)$ covariant
 valuations leads to an algebraic characterization of all
 $\mathcal{G}(\Z^2)$ covariant
 valuations.

\begin{lemma}
\label{simplef2algebraic-properties}
For any simple $\mathcal{G}(\Z^2)$ covariant valuation $Z:\mathcal{P}(\Z^2)\to \mathcal{M}(\R^2)$,
$Z(T)=f_2$ satisfies the following properties: For $x,y\in\R$, we have
\begin{align}
\label{f2simple1}
f_2(-x+y,-x)&= e^{-x}f_2(x,y);\\
\label{f2simple2}
f_2(x,y)+e^{x+y}f_2(-x,-y)&=f_2(x,x+y)+f_2(x+y,y);\\
\label{f2simple3}
f_2(x,y)&= f_2(y,x).
\end{align}
\end{lemma}
\proof The $\mathcal{G}(\Z^2)$ covariance yields
\eqref{f2simple1} because
$T-e_1=\left[
\begin{array}{cc}
-1&-1\\
1&0
\end{array}\right] T$.

For \eqref{f2simple2},
we observe that each diagonal of the square
$[0,1]^2$ cuts the square into two empty triangles. In particular,
$[0,1]^2$
can be written as $T\cup (e_1+e_2-T)$ on the one hand, and the union of the triangles
 $\left[
\begin{array}{cc}
1&1\\
0&1
\end{array}\right] T$
and
$\left[
\begin{array}{cc}
1&0\\
1&1
\end{array}\right] T$ on the other hand.
 We conclude \eqref{f2simple2} from evaluating $Z([0,1]^2)$ in the two ways; and
using the simpleness of $Z$ and the $\mathcal{G}(\Z^2)$ covariance of $Z$.

Finally, \eqref{f2simple3} follows from $T=\left[
\begin{array}{cc}
0&1\\
1&0
\end{array}\right] T$ and
 the $\mathcal{G}(\Z^2)$ covariance of $Z$.
\endproof

\begin{prop}
\label{simplef2algebraic}
For any $f_2\in\mathcal{M}(\R^2)$ satisfying
the properties \eqref{f2simple1}, \eqref{f2simple2} and \eqref{f2simple3}
in Lemma~\ref{simplef2algebraic-properties},
there exists a unique simple $\mathcal{G}(\Z^2)$ covariant
valuation $Z:\mathcal{P}(\Z^2)\to \mathcal{M}(\R^2)$ such that
$Z(T)=f_2$.
\end{prop}
\proof We write $H$  to denote the subgroup of stabilizers of $T$ of
$\mathcal{G}(\Z^2)$, and let $H_+\subset H$ be the subgroup of orientation preserving transformations.

In accordance with Lemma~\ref{StabilizerGenertors}, $H_+$
 has three elements corresponding to the vertices of
$T$; namely, the identity, $\Xi_+x=\Phi_+x+e_1$
for $\Phi_+=\left[
\begin{array}{cc}
-1&-1\\
1&0
\end{array}\right]$
where $\Xi_+T=T$, and $\Xi_+^{-1}x=\Phi_+^{-1}x+e_2$ where
 $\Phi_+^{-1}=\left[
\begin{array}{cc}
0&1\\
-1&-1
\end{array}\right]$. It follows from Lemma~\ref{StabilizerGenertors} (ii) that $H$ is generated by
$\Xi_+$ and $\Phi_{-}=\left[
\begin{array}{cc}
0&1\\
1&0
\end{array}\right]$. Since the property \eqref{f2simple1} of $f_2$
 corresponds to the $\Xi_+$ invariance of $T$, and
 the property \eqref{f2simple3} of $f_2$
 corresponds to the $\Phi_{-}$ invariance of $T$, it follows from
Lemma~\ref{well-defined}  that we may define
$Z_2(N)$ for any empty triangle $N$ by choosing a $\Xi\in \mathcal{G}(\Z^2)$ with $N=\Xi T$ and setting
\begin{equation}
\label{Z2EmptyTriangledef}
Z_2(N)=\Xi\cdot f_2.
\end{equation}
In addition, the $Z_2$ defined as in \eqref{Z2EmptyTriangledef} is $\mathcal{G}(\Z^2)$  covariant on empty lattice triangles
according to Lemma~\ref{well-defined}.

Next we call a parallelogram $P\in\mathcal{P}(\Z^2)$ an empty parallelogram if contains no other lattice points
then its vertices. In this case, $P=\Xi[0,1]^2$ for a $\Xi\in \mathcal{G}(\Z^2)$.

For an empty parallelogram $P\in\mathcal{P}(\Z^2)$, one of the diagonals of $P$ cuts $P$ into the
 two empty triangles, let them be $T_1$ and $T_2$, and
let $T_3$ and $T_4$ be the two empty triangles in $P$ determined by the other diagonal of $P$.
We claim that
\begin{equation}
\label{Z2EmptyParallelogram}
Z_2(T_1)+Z_2(T_2)=Z_2(T_3)+Z_2(T_4).
\end{equation}
Let $v$ be the vertex of $T_1$ that is not a vertex of $T_2$, and let $\Xi\in \mathcal{G}(\Z^2)$ such that
$P=\Xi[0,1]^2$ and $\Xi o=v$.  For $\widetilde{\Xi}_1=I_2$ and $\widetilde \Xi_2x=-x+e_1+e_2$,
the diagonal $[e_1,e_2]$ of $[0,1]^2$ divides the square into $\widetilde{\Xi}_1T$
and $\widetilde{\Xi}_2T$, and $T_i=\Xi\,\widetilde{\Xi}_iT$ for $i=1,2$.
Possibly after interchanging $T_3$ and $T_4$, we may assume that
$T_i=\Xi\,\widetilde{\Xi}_iT$ for $i=3,4$ where
 $\widetilde{\Xi}_3=\left[
\begin{array}{cc}
1&1\\
0&1
\end{array}\right]$
and
$\widetilde{\Xi}_4=\left[
\begin{array}{cc}
1&0\\
1&1
\end{array}\right]$, and the diagonal $[o,e_1+e_2]$ of $[0,1]^2$ divides the square into
$\widetilde{\Xi}_3T$
and $\widetilde{\Xi}_4T$. Since \eqref{f2simple2} states that
$$
\widetilde{\Xi}_1\cdot f_2+\widetilde{\Xi}_2\cdot f_2=\widetilde{\Xi}_3\cdot f_2+\widetilde{\Xi}_4\cdot f_2,
$$
we conclude \eqref{Z2EmptyParallelogram} from the definition and $\mathcal{G}(\Z^2)$ covariance of $Z_2$ for empty triangles.

To define $Z_2(Q)$ for a lattice polygon $Q\in\mathcal{P}(\Z^2)$, to any
triangulation $\mathcal{T}=\{E_1,\ldots,E_k\}$ into empty triangles, we assign
$$
\xi(\mathcal{T})=\sum_{i=1}^kZ_2(E_i).
$$
We claim that if $\mathcal{T}'=\{E'_1,\ldots,E'_k\}$ is another triangulation of $Q$ into empty triangles, then
\begin{equation}
\label{PolygonTriangulation}
\xi(\mathcal{T})=\xi(\mathcal{T}').
\end{equation}
According to Lawson \cite{Law72}
(see J. de Loera, J. Rambau, F. Santos \cite{LRS10}, Section 3.4.1), the two triangulations $\mathcal{T}$
and $\mathcal{T}'$ can be connected by a series of triangulations  $\mathcal{T}_0,\ldots,\mathcal{T}_m$
of $Q$ into empty triangles such that $\mathcal{T}=\mathcal{T}_0$, $\mathcal{T}'=\mathcal{T}_m$,
and $\mathcal{T}_{i+1}$ is obtained by a diagonal flip from $\mathcal{T}_i$. Here a diagonal flip means that
we find two empty triangles $T_1$ and $T_2$ in  $\mathcal{T}_i$ whose union is an empty parallelogram $P$,
and hence $T_1\cap T_2$ is a diagonal of $P$,
and we replace $T_1$ and $T_2$ by the two empty triangles $T_3$ and $T_4$ obtained by cutting $P$ into two by the
other diagonal of $P$. We deduce from \eqref{Z2EmptyParallelogram} that
$\xi(\mathcal{T}_{i+1})=\xi(\mathcal{T}_i)$, yielding \eqref{PolygonTriangulation}.

It follows from \eqref{PolygonTriangulation} that
 for any lattice polygon $Q\in\mathcal{P}(\Z^2)$, we can define $Z_2(Q)$ by subdividing $Q$ into
empty lattice triangles $E_1,\ldots,E_k$, and setting
\begin{equation}
\label{Z2Polygondef}
Z_2(Q)=\sum_{i=1}^kZ_2(E_i).
\end{equation}
If $P\in\mathcal{P}(\Z^2)$ is a point or a segment, then we define $Z_2(P)=0$. Since
$Z_2$ has been known to be $\mathcal{G}(\Z^2)$ contravariant on empty triangles, $Z_2$ is  $\mathcal{G}(\Z^2)$ contravariant
on $\mathcal{P}(\Z^2)$.

Therefore all left to prove is that $Z_2$ is a simple valuation. To verify this, we only need to check the case
when $P,Q\in\mathcal{P}(\Z^2)$ are lattice polygons such that $P\cup Q$ is convex and
$P\cap Q$ is also a lattice polygon. First we subdivide $P\cap Q$ into empty triangles, which triangulation we extend
into a triangulation $\mathcal{T}_P$  of $P$ and a triangulation $\mathcal{T}_Q$  of $Q$ into empty triangles.
It follows that $\mathcal{T}_P\cap\mathcal{T}_Q$ is the original triangulation of $P\cap Q$,
and $\mathcal{T}_P\cup\mathcal{T}_Q$ is a triangulation of $P\cup Q$. Therefore,
\eqref{Z2Polygondef} yields that $Z_2(P\cup Q)+Z_2(P\cap Q)=Z_2(P)+Z_2(Q)$.
\endproof

We note that for any $a,b\in \N$, the function
\begin{equation}
\label{x-ydivideanalytic}
(x,y)\mapsto \frac1{x-y}\left(\frac{e^{y}-1}{y}\cdot x^ay^b-\frac{e^{x}-1}{x}\cdot x^by^a\right)
\end{equation}
is analytic.

\begin{example}
\label{analytic-example}
Since the analytic function ({\it cf.} \eqref{x-ydivideanalytic})
\begin{align*}
f_2(x,y)=&  \frac1{x-y}\cdot \left(\frac{e^{y}-1}{y}\cdot(x^4-4x^3y+x^2y^2+6xy^3-3y^4)\right.-\\
& \left.-\frac{e^{x}-1}{x}\cdot (y^4-4xy^3+x^2y^2+6x^3y-3x^4)\right)
\end{align*}
satisfies \eqref{f2simple1}, \eqref{f2simple2} and \eqref{f2simple3}, it
follows from Proposition~\ref{simplef2algebraic} that there exists a
$\mathcal{G}(\Z^2)$  covariant simple
valuation $Z:\mathcal{P}(\Z^2)\to \mathcal{M}(\R^2)$ satisfying
$Z(T)=f_2$ that is
different from the Laplace transform (see \eqref{LaplaceTformulaxt}), and $Z(P)$ is analytic for any lattice polytope $P$.
In the language of  Freyer, Ludwig, Rubey \cite{FMR}, $Z=\sum_{i>0}L^{i+2}_i$ where $L^r_i$ is the $i$-homogeneous Ehrhart tensor on lattice polygons with rank $r$.
\end{example}

\section{$Z(\{o\})=f_0$ is essentially constant}
\label{secSL2Zergodic}

The action of a group G on a space $S$ equipped with a measure $\mu$ is called {\it ergodic} if for any $G$ invariant measurable $X\subset S$, either $\mu(X)=0$ or $\mu(S\backslash X)=0$.
The following well-known statement is Example~2.2.9 in Zimmer \cite{Zim84}.

\begin{prop}
\label{SL2Zergodic}
The action of ${\rm SL}(2,\Z)$ on $\R^2$ is ergodic.
\end{prop}

Proposition~\ref{SL2Zergodic} directly yields Proposition~\ref{fSL2Zinvariant}, that here we recall as
Corollary~\ref{fSL2Zinvariant0},
stating that any ${\rm SL}(2,\Z)$ invariant measurable function is essentially constant
(see Proposition~2.1.11 in Zimmer \cite{Zim84}).

\begin{coro}
\label{fSL2Zinvariant0}
If $f\in \mathcal{M}(\R^2)$ is invariant under ${\rm SL}(2,\Z)$, then there exists a constant $c\in\R$
such that $f(x)=c$ for almost everywhere $x\in\R^2$.
\end{coro}
\proof Let $M=\sup\{t:|\{f\geq t\}|> 0\}\in(-\infty,\infty]$ and $m=\inf\{t:|\{f\leq t\}|> 0\}\in[-\infty,\infty)$,
and hence $m\leq M$. If $m=M$, then $c=m=M$ works as $|\{f>c\}|= 0$ and $|\{f<c\}|= 0$ (here we did not even use
the ${\rm SL}(2,\Z)$ invariance of $f$).

If $m<M$, then choose $c\in(m,M)$. Since exactly one of the pairwise disjoint ${\rm SL}(2,\Z)$ invariant measurable sets
$\{f>c\}$, $\{f=c\}$ and $\{f<c\}$ has non-zero measure, and the union of the first two and the union of the last two have positive measure, we have $|\{f>c\}|=|\{f<c\}|=0$.
\endproof

Combining Lemma~\ref{f0f1f2properties} and Corollary~\ref{fSL2Zinvariant0} yields the characterization of $f_0=Z(\{o\})$.

\begin{prop}
\label{f0characterization}
For any translatively exponential  and ${\rm GL}(2,\Z)$ covariant
valuations $Z:\mathcal{P}(\Z^2)\to \mathcal{M}(\R^2)$, $f_0=Z(\{o\})$ is ${\rm GL}(2,\Z)$ invariant, and hence almost everywhere constant.
\end{prop}

\section{Characterization of $Z([o,e_1])=f_1$, and the proof of Theorem~\ref{ExpValGL2Z}}
\label{secf1}

According to Lemma~\ref{f0f1f2properties} and Proposition~\ref{lowdim},
for an $f_1\in\mathcal{M}(\R^2)$, there exists a  $\mathcal{G}(\Z^2)$ covariant valuation $Z:\mathcal{P}(\Z^2)\to\mathcal{M}(\R^2)$ satisfying
$Z([o,e_1])=f_1$ if and only if
for any $(x,y)\in\R^2$, we have
\begin{eqnarray}
\label{f1shift0}
f_1(x,y)&=& e^{x}f_1(-x,-y);\\
\label{f11110}
f_1(x,y)&=& f_1(x,x+y),\\
\label{f1extra0}
f_1(x,y)&=& f_1(x,-y).
\end{eqnarray}

\begin{prop} \mbox{ }
\label{f1geometric}
 For $\Omega_{1}=\{(0,y):y\geq 0\}\cup \{(x,y):0\leq y\leq \frac12\,x\}$,
and for any measurable $\tilde{f}_1:\Omega_{1}\to\R$, there exists a unique
$f_1\in\mathcal{M}(\R^2)$ satisfying \eqref{f1shift0}, \eqref{f11110}, \eqref{f1extra0} and
 $f_1|_{\Omega_{1}}=\tilde{f}_1$.
\end{prop}
\proof The key observation is that \eqref{f11110} yields that
\begin{equation}
\label{f1111k}
f_1(x,y)=f_1(x,kx+y)\mbox{ \ for $k\in\Z$ and $x,y\in\R$.}
\end{equation}

Now knowing $f_1$ on $\{(x,y):0\leq y\leq \frac12\,x\}$
determines
$f_1(x,y)$ for $x>0$ and $|y|\leq\frac12\,x$ via \eqref{f1extra0}.
We observe that
\begin{equation}
\label{f1GLwell-defined}
f_1\left(x,-\mbox{$\frac12\,$}x\right)=f_1\left(x,\mbox{$\frac12\,$}x\right)
\mbox{ \ for $x>0$.}
\end{equation}
Next  knowing $f_1$ on $\{(x,y):|y|\leq \frac12\,x\}$
determines
$f_1(x,y)$ for $x>0$ and $y\in\R$ via \eqref{f1111k} where
\eqref{f1GLwell-defined} ensures that $f_1$ is well-defined. In addition, knowing
$f_1$ on $\{(0,y):\,y\geq 0\}$ determines $f_1(0,y)$ for $y\in\R$ as \eqref{f1shift0}
implies $f_1(0,-y)=f_1(0,y)$. Then,
knowing
$f_1$ on $\{(x,y):x> 0\}$ determines $f_1(x,y)$ for $x<0$ and $y\in\R$ via
\eqref{f1shift0}. Finally, the $f_1$ constructed this way satisfies
\eqref{f1shift0}, \eqref{f11110} and \eqref{f1extra0}.
\endproof

\proof[Proof of Theorem~\ref{ExpValGL2Z}.]
Combining Proposition~\ref{lowdim}, Corollary~\ref{ReductionToSimple}, Proposition~\ref{f0characterization}
and Proposition~\ref{f1geometric} implies
Theorem~\ref{ExpValGL2Z}.
\endproof

\section{Towards a Geometric Characterization of $Z(T)=f_2$ for simple valuations}
\label{secf2step1}

According to Lemma~\ref{simplef2algebraic-properties} and Proposition~\ref{simplef2algebraic} in Section~\ref{secsimplef2},  if  $Z:\mathcal{P}(\Z^2)\to \mathcal{M}(\R^2)$ is a $\mathcal{G}(\Z^2)$ covariant simple
valuation, then
$Z(T)=f_2$ is characterized by the following properties (cf. \eqref{f2simple1}, \eqref{f2simple2} and \eqref{f2simple3}
in Lemma~\ref{simplef2algebraic-properties}):
\begin{align}
\label{f220}
f_2(-x+y,-x)&= e^{-x}f_2(x,y);\\
\label{f230}
f_2(x,y)+e^{x+y}f_2(-x,-y)&=f_2(x,x+y)+f_2(x+y,y);\\
\label{f210}
f_2(x,y)&= f_2(y,x).
\end{align}

Next we provide a partially geometric characterization of functions satisfying \eqref{f220}, \eqref{f230} and \eqref{f210}. As the first step, we concentrate only on the properties \eqref{f220} and \eqref{f210}. Let
$$
\widetilde{\Omega}=\{(x,y)\in\R^2:\,y\geq x\geq 0\}.
$$
We consider the positive hull ${\rm pos}\{u,v\}=\{\alpha u+\beta v:\alpha,\beta\geq 0\}$ of $u,v\in\R^2$.

\begin{lemma}
\label{iandii}
Given any measurable function $\tilde{f}:\,\widetilde{\Omega}\to\R$, there exists a unique
 measurable $f_2:\,\R^2\to\R$ satisfying \eqref{f220} and  \eqref{f210} extending $\tilde{f}$.
\end{lemma}
\proof

We consider the $\Phi\in{\rm SL}(n,\Z)$ defined by
$$
\Phi(x,y)=(-x+y,-x),
$$
and hence $\Phi^{-1}(x,y)=(-y,x-y)$, 
and \eqref{f220} is equivalent to 
\begin{equation}
\label{f220Phi}
\begin{array}{rcl}
f_2(\Phi(x,y))&=&e^{-x}f_2(x,y)\\[1ex]
e^{-y}f_2(x,y)&=&f_2(\Phi^{-1}(x,y))=f_2(-y,x-y).
\end{array}
\end{equation}
We observe that $\Phi$ maps the quadrant $[0,\infty)^2={\rm pos}\{e_1,e_2\}$ onto ${\rm pos}\{e_1,-e_1-e_2\}$, while $\Phi^{-1}=\Phi^2$
maps $[0,\infty)^2$ onto ${\rm pos}\{e_2,-e_1-e_2\}$.

Given a measurable function $\tilde{f}:\,\widetilde{\Omega}\to\R$, we extend it to an
$f$ on $[0,\infty)^2$ by setting $f(x,y)=\tilde{f}(y,x)$ if $0\leq y\leq x$. This function
satisfies that $f(\Phi(x,y))=e^{-x}f(x,y)$ if $(x,y),\Phi(x,y)\in[0,\infty)^2$ because $x=0$ in this case. Using the properties of $\Phi$, we can extend $f$ in a unique way to a function $f_2$ on $\R^2$ satisfying
\eqref{f220Phi}, and hence \eqref{f220}.

Finally, we need to check whether $f_2$ satisfies \eqref{f210}. It does hold on $[0,\infty)^2$ by definition,
and the $\Phi[0,\infty)^2$
and $\Phi^{-1}[0,\infty)^2$ are images of each other by the linear transformation
$(x,y)\mapsto (y,x)$. Therefore, all we need to check is that
$f_2(x,y)=f_2(y,x)$ for
 $(x,y)\in \Phi^{-1}[0,\infty)^2$ with $x\neq y$, and hence $x<0$ and $y>x$.
We have $(x,y)=\Phi^{-1}(y-x,-x)$ and $(y,x)=\Phi(-x,y-x)$, and hence \eqref{f220Phi} implies
\begin{align*}
f_2(x,y)&=e^xf(y-x,-x)\\
f_2(y,x)&=e^xf(-x,y-x)
\end{align*}
where $f(y-x,-x)=f(-x,y-x)$ by the construction of $f$.
\endproof

Based on Lemma~\ref{iandii}, we plan to include the condition  \eqref{f230} in the characterization, and hence we rewrite \eqref{f230} into forms suiting better the domains in Lemma~\ref{iandii}.

\begin{lemma}
\label{f23uplemma}
 Assuming that $f_2:\,\R^2\to\R$ satisfies \eqref{f220} and \eqref{f210}, we have $f_2$ also satisfies
\eqref{f230} if and only if for any $x,y\in \R$, we have
\begin{equation}
\label{f23up}
f_2(x,y)+e^{x}f_2(y-x,y)=f_2(x,x+y)+f_2(y,x+y).
\end{equation}
\end{lemma}
\proof
 We observe that using the substitution
$a=-x+y$ and $-x=b$,  \eqref{f220} is equivalent with
$e^{-b}f_2(a,b)=f_2(-b,a-b)$, and in turn \eqref{f220} is equivalent with
 $e^{y}f_2(-x,-y)=f_2(y,y-x)$ for any $x,y\in\R$. Therefore, under the condition
of symmetry; namely, that $f$ satisfies \eqref{f210}, the condition \eqref{f220} is equivalent with
$e^{x+y}f_2(-x,-y)=e^xf_2(y-x,y)$ for any $x,y\in\R$, and then \eqref{f230} is equivalent to \eqref{f23up}.
\endproof

\section{The proof of Theorem~\ref{simpleGLcovariant}}
\label{secf2full}

Concerning \eqref{f23up}, we observe that if $(x,y)\in\widetilde{\Omega}$
where
$$
\widetilde{\Omega}=\{(x,y)\in\R^2:\,y\geq x\geq 0\},
$$
then
$(y-x,y)$, $(x,x+y)$ and $(y,x+y)$ all lie in $\widetilde{\Omega}$. Therefore,
it makes sense to consider \eqref{f23up} as a property of a function on $\widetilde{\Omega}$.

\begin{prop}
\label{f2GeometricAlgebraicGL}
 Given any measurable function $\tilde{f}:\,\widetilde{\Omega}\to\R$ satisfying
\eqref{f23up}, there exists a unique
 measurable $f_2:\,\R^2\to\R$ satisfying \eqref{f220}, \eqref{f230} and  \eqref{f210} extending $\tilde{f}$.
\end{prop}
\proof Using Lemma~\ref{iandii}, first we extend $\tilde{f}$ to the unique measurable function $f_2$ on
$\R^2$ satisfying \eqref{f220} and  \eqref{f210}, and we
 use \eqref{f220}
in both forms
\begin{equation}
\label{fr220both}
e^{-a}f_2(a,b)=f_2(-a+b,-a)\mbox{ and }e^{-b}f_2(a,b)=f_2(-b,a-b).
\end{equation}
We note that for $f_2$, the formulas
$f_2(a,b)=f_2(b,a)$ and \eqref{fr220both} yield
\begin{align}
\nonumber
f_2(x,y)+e^{x}f_2(y-x,y)&=f_2(y,x)+e^{x}e^{y-x}f_2(x,x-y)\\
\label{f23upLeftSwitch}
&=f_2(y,x)+e^{y}f_2(x-y,x).
\end{align}

Now we
verify that  $f_2$ satisfies \eqref{f23up} for any $x,y\in\R$.

\noindent{\bf Case 1} $x,y\geq 0$\\
 If $y\geq x$, then
$(x,y)\in\widetilde{\Omega}$, and \eqref{f23up} holds. If $0\leq y<x$, then combining  \eqref{f23up} for $(y,x)\in\widetilde{\Omega}$ and 
\eqref{f23upLeftSwitch}
implies \eqref{f23up} again.\\

\noindent{\bf Case 2} $x\leq 0$ and $y\geq 0$\\
Now $(-x,y)\in[0,\infty)^2$; therefore, Case 1 yields that
$$
f_2(-x,y)+e^{-x}f_2(y+x,y)=f_2(-x,y-x)+f_2(y,y-x).
$$
We multiply through by $e^x$, and observe that  \eqref{fr220both} and $f_2(a,b)=f_2(b,a)$ imply
$$
e^{x}f_2(-x,y)+f_2(y+x,y)=f_2(x+y,x)+f_2(y+x,y)=f_2(x,x+y)+f_2(y,y+x),
$$
 which is the right hand side of \eqref{f23up}. On the other hand,
again \eqref{fr220both} and $f_2(a,b)=f_2(b,a)$ yield
$$
e^xf_2(-x,y-x)+e^xf_2(y,y-x)=f_2(y,x)+e^xf_2(y,y-x)=f_2(x,y)+e^xf_2(y-x,y),
$$
which is the left hand side of \eqref{f23up}. In turn, we conclude  \eqref{f23up}.\\

\noindent{\bf Case 3} $x\geq 0$ and $y\leq 0$\\
Now Case 2 yields that
$$
f_2(y,x)+e^{y}f_2(x-y,x)=f_2(y,y+x)+f_2(x,y+x),
$$
and hence \eqref{f23upLeftSwitch} implies \eqref{f23up}.\\

\noindent{\bf Case 4} $x,y< 0$\\
As $-x> 0$ and $y<0$,  Case 3 yields that
$$
f_2(-x,y)+e^{-x}f_2(y+x,y)=f_2(-x,y-x)+f_2(y,y-x).
$$
Therefore, the argument in Case 2 yields \eqref{f23up}.
\endproof

According to Proposition~\ref{f2GeometricAlgebraicGL}, the task to construct  
 a measurable function $f_2:\,\R^2\to\R$ satisfying \eqref{f220}, \eqref{f230} and  \eqref{f210}
- that is needed for Theorem~\ref{simpleGLcovariant} by Lemma~\ref{simplef2algebraic-properties} and Proposition~\ref{simplef2algebraic} -
is equivalent to construct a measurable function
$\tilde{f}:\,\widetilde{\Omega}\to\R$ satisfying
\begin{equation}
\label{f23uptilde}
\tilde{f}(x,y)+e^{x}\tilde{f}(y-x,y)=\tilde{f}(x,x+y)+\tilde{f}(y,x+y).
\end{equation}
When we write $\varphi(t)=\frac{e^t-1}{t}$ in the formulas below, we mean the positive analytic function
on $\R$ with $\varphi(0)=1$; namely, $\varphi(t)=\sum_{n=0}^\infty\frac{t^n}{(n+1)!}$, and hence the reciprocal of $\varphi$ is the series defining the Bernoulli numbers $B_0, B_1,\ldots$ by the formula (cf. Zagier \cite{Zag08})
$$
\frac{t}{e^t-1}=\sum_{n=0}^\infty\frac{B_n}{n!}\cdot t^n.
$$
The following statement follows by direct calculations.

\begin{lemma} 
\label{tildefrhorelation}
A natural bijection  between 
measurable functions
$\tilde{f}:\,\widetilde{\Omega}\to\R$ satisfying \eqref{f23uptilde}
and $\varrho:\,[0,\infty)^2\to\R$
satisfying \eqref{rhoeq0}
is induced by \eqref{rhotildef}
and \eqref{tildefrho} as follows:
\begin{description}
\item{(i)}
 If a measurable function $\tilde{f}:\,\widetilde{\Omega}\to\R$ satisfies
\eqref{f23uptilde}, then the function $\varrho:\,[0,\infty)^2\to\R$ defined by
\begin{equation}
\label{rhotildef}
\varrho(x,y)=\frac{x}{e^x-1}\cdot
\frac{x+y}{e^{x+y}-1}\cdot \left(\tilde{f}(x,x+y)+e^x\cdot \tilde{f}(y,x+y)\right)
\end{equation}
  satisfies
\begin{equation}
\label{rhoeq0}
(2x+y) \varrho(x,y)=(x+y)\varrho(x,x+y)+x \varrho(x+y,x),
\end{equation}
$\varrho(0,0)=2\tilde{f}(0,0)$, and for $(x,y)\in \widetilde{\Omega}\backslash(0,0)$, the formula
\begin{equation}
\label{tildefrho}
\tilde{f}(x,y)=\frac{e^x}{y}\cdot \frac{e^{y-x}-1}{y-x}\cdot \varrho(y-x,x)-
\frac{1}y\cdot \frac{e^{x}-1}{x}\cdot \varrho(x,y-x).
\end{equation}
\item{(ii)}
If a function $\varrho:\,[0,\infty)^2\to\R$ satisfies \eqref{rhoeq0}, then
the function $\tilde{f}:\,\widetilde{\Omega}\to\R$ defined by \eqref{tildefrho}
 satisfies
\eqref{f23uptilde}.
\end{description}
\end{lemma}
\noindent{\bf Remark.}
The case when the original simple valuation is the positive  Laplace transform is equivalent with the property that
 $\varrho$ is the constant $1$ function in \eqref{rhoeq0} (cf. \eqref{LaplaceTformulaxt}).
 \proof
Substituting the formula \eqref{tildefrho} for $\tilde{f}$ into the left side of \eqref{f23uptilde} results in the expression $\frac{(e^x - 1)(e^y - 1)}{xy} \varrho(x,y-x)$, and into 
 the right side of \eqref{f23uptilde} results in the expression
  $$
 \frac{(e^x - 1)(e^y - 1)}{xy(x+y)}(y\varrho(x,y) + x \varrho (y,x)).
 $$ 
 It follows that
 $$
 y\varrho(x,y) + x \varrho (y,x) = (x+y)\varrho(x,y-x),
 $$
 which is equivalent to \eqref{rhoeq0} by
  replacing $y$ by $y-x$. This proves (ii), and (i) can be verified analogously.  
 \endproof

To provide a characterization of the measurable function $\varrho:\,[0,\infty)^2\to\R$
in Lemma~\ref{tildefrhorelation}, we consider the golden ratio
$$
\tau=\frac{\sqrt{5}+1}2,
$$
which satisfies the properties
\begin{equation}
\label{golden-ratio-rho}
\tau^2=1+\tau\mbox{ and } \frac1{\tau}=\tau-1\mbox{ and } \frac{1+\tau}{\tau}=\tau.
\end{equation}
In addition, recall that
\begin{equation}
\label{Omega2def}
\widetilde{\Omega}_{2}=\{(\tau s,s):1\leq s<\tau\}\cup \{(x,y):0\leq x\leq y\}.
\end{equation}

\begin{lemma}
\label{Fibonacci-rho}
If $\tilde{\varrho}$ is any (or Lebesgue measurable, or Borel measurable) function on $\widetilde{\Omega}_2$ into a real vectorspace, then there exists a unique (Lebesgue measurable, or Borel measurable) function $\varrho$ on
$[0,\infty)^2$ satisfying
\begin{equation}
\label{rhoformula0}
(2x+y) \varrho(x,y)=(x+y)\varrho(x,x+y)+x \varrho(x+y,x),
\end{equation}
whose restriction to $\widetilde{\Omega}_2$ is $\tilde{\varrho}$.
\end{lemma}
\proof  If $x=0$ and $y\geq  0$, then 
both the left hand side and the right hand side of \eqref{rhoformula0} is
$y\varrho(0,y)$, thus \eqref{rhoformula0} provides no restriction on $\varrho$.
 If $x>0$ and $y=0$, then \eqref{rhoformula0} reads as $2x\varrho(x,0)=2x\varrho(x,x)$; therefore,
$\varrho(x,0)=\varrho(x,x)$. In particular, we may assume that $x,y>0$ for
$(x,y)\in\widetilde{\Omega}_{2}$ in \eqref{rhoformula0}. Since then $(x+y,x)\not\in \widetilde{\Omega}_{2}$,
\eqref{rhoformula0} provides no restriction on $\tilde{\varrho}$ on $\widetilde{\Omega}_{2}$.
However, we still need to define $\varrho$ on $\Omega_*=\{(x,y):0<y<x\}$.

We recall the Fibonacci sequence, which we index as
$$
F_0=0, \mbox{ and } F_1=1, \mbox{ and $F_{n+2}=F_n+F_{n+1}$ for $n\geq 0$.}
$$
For $n\geq 0$, setting $\frac1{0}=\infty$, the Fibonacci sequence (cf. Burton \cite{Bur89}) satisfies
\begin{eqnarray*}
\frac{F_{n+1}}{F_{n}}&>&\frac{F_{n+3}}{F_{n+2}}>\tau\mbox{ \ if $n\geq 0$ even;}\\
\frac{F_{n+1}}{F_{n}}&<&\frac{F_{n+3}}{F_{n+2}}<\tau\mbox{ \ if $n\geq 1$ odd;}\\
\lim_{n\to\infty}\frac{F_{n+1}}{F_{n}}&=&\tau.
\end{eqnarray*}
In particular, for $n\geq 0$,  we can define
\begin{eqnarray*}
\Omega_n&=&\left\{(x,y):\, \frac{F_{n+1}}{F_{n}}>\frac{x}{y}\geq \frac{F_{n+3}}{F_{n+2}}\right\}
\mbox{ \ if $n\geq 0$ is even};\\
\Omega_n&=&\left\{(x,y):\, \frac{F_{n+1}}{F_{n}}< \frac{x}{y}\leq \frac{F_{n+3}}{F_{n+2}}\right\}
\mbox{ \ if $n\geq 1$ is odd.}
\end{eqnarray*}
 We observe that $\{\Omega_n\}_{n\geq 0}$ is a partition of $\Omega_*\backslash \ell$
for the open half line $\ell= \{(a,b):\,a=\tau b\mbox{ and }b>0\}$
 by the properties above of he Fibonacci numbers.

We note that for the points of $\Omega_*$  and for $n\geq 0$, we have
\begin{eqnarray}
\label{rho-1}
\mbox{if }(x,y)\in\Omega_*,&\mbox{ then }&(x,x+y)\in\widetilde{\Omega}_{2};\\
\label{rho-2}
(x,y)\in\Omega_n&\mbox{ if and only if }&(x+y,x)\in\Omega_{n+1}\\
\label{rho-yxtau}
x,y>0\mbox{ and }\frac{x}{y}=\tau &\mbox{ if and only if }& y>0\mbox{ and }\frac{x+y}{x}=\tau.
\end{eqnarray}
We observe that for any $x,y$, $(x+y,x)\in\Omega_{0}=\{(a,b):a\geq 2b>0\}$ if and only
if $(x,y)\in\{(a,b):0<a\leq b\}\subset \widetilde{\Omega}_{2}$; therefore,
we can define $\varrho$ on $\Omega_0$ {\it via} \eqref{rhoformula0}. It can be done in a unique way because
if $(x,y)\in\Omega_0$, then $(x+y,x)\in\Omega_{1}$ by \eqref{rho-2}.
Using \eqref{rho-1} and \eqref{rho-2}, we define $\varrho$ on any $\Omega_n$, $n\geq 0$, by induction on $n\geq 0$. This extension of $\varrho$ onto $\Omega_*\backslash \ell$ is also unique
 according to \eqref{rho-1} and \eqref{rho-2}.
Therefore all we need to do is to define $\varrho$ on the open half line
$\ell$.

For $m\in\Z$, let $s_m=\{(a,b):\,a=\tau b\mbox{ and }\tau^m\leq b<\tau^{m+1}\}\subset\ell$, and hence
$\{s_m\}_{m\in\Z}$ provides a partition of $\ell$, and $\ell\cap \widetilde{\Omega}_{2}=s_0$.
Since $(x,y)\in s_m$ if and only if $(x+y,x)\in s_{m+1}$ according to \eqref{rho-yxtau},
we can define $\varrho$ on $s_n$ for $n\geq 0$ by induction on $n$ {\it via}
\eqref{rhoformula0}, and similarly we can define $\varrho$ on $s_{-n}$ for $n\geq 0$ by induction on $n$
again {\it via}
\eqref{rhoformula0}, completing the proof of Lemma~\ref{Fibonacci-rho}.
\endproof

\proof[Proof of Theorem~\ref{simpleGLcovariant}.]
Combining Lemma~\ref{simplef2algebraic-properties}, Proposition~\ref{simplef2algebraic}, Proposition~\ref{f2GeometricAlgebraicGL}, Lemma~\ref{tildefrhorelation} and Lemma~\ref{Fibonacci-rho}
yields Theorem~\ref{simpleGLcovariant}.
\endproof

\bigskip

\noindent{\bf Acknowledgements.} We are grateful to Vitaly Bergelson (Ohio State University), Lewis Bowen (University of Texas at Austin), Monika Ludwig  (TU Wien) and \'Arp\'ad T\'oth (ELTE) for enlightening discussions. We are grateful for the referee for signicantly improving the presentation and the focus of the paper.

K\'aroly J. B\"or\"oczky, HUN-REN Alfr\'ed R\'enyi Institute of Mathematics, boroczky.karoly.j@renyi.hu\\

M\'aty\'as Domokos, HUN-REN Alfr\'ed R\'enyi Institute of Mathematics, domokos.matyas@renyi.hu\\

Ansgar Freyer, FU Berlin, Fachbereich Mathematik und Informatik, Arnimallee 2, 14195 Berlin, a.freyer@fu-berlin.de\\

Christoph Haberl, Vienna University of Technology, Institute of Discrete Mathematics and Geometry, Wiedner Hauptstraße 8-10/104, 1040 Vienna, Austria, christoph.haberl@tuwien.ac.at\\

Gergely Harcos, HUN-REN Alfr\'ed R\'enyi Institute of Mathematics, harcos.gergely@renyi.hu\\

Jin Li, Department of Mathematics, and Newtouch Center for Mathematics, Shanghai University, li.jin.math@outlook.com 

\end{document}